\documentclass[titlepage,11pt]{article}
\usepackage{latexsym}
\usepackage{amsfonts}
\usepackage{amsmath}
\usepackage{comment}
\newcommand\blackslug{\hbox{\hskip 1pt \vrule width 4pt height 8pt depth 1.5pt
        \hskip 1pt}}
\newcommand\bbox{\hfill \quad \blackslug \bigbreak}

\def\ll{,\ldots,}
\newcommand\vare{\varepsilon}

\title{Pure pairs. VI. Excluding an ordered tree}
\author{Alex Scott
\\
Mathematical Institute, University of Oxford, Oxford OX2 6GG, UK
\\
\\
Paul Seymour\thanks{Supported by AFOSR grant
A9550-19-1-0187 and NSF grant  DMS-1800053.}\\
Princeton University, Princeton, NJ 08544
\\
\\
Sophie Spirkl\thanks{This material is based upon work supported by the National Science
Foundation under Award No. DMS-1802201.}\\
Princeton University, Princeton, NJ 08544}

\date{March 10, 2020; revised \today}

\newtheorem{thm}{}[section]

\newcommand{\Proof}{\noindent{\bf Proof.}\ \ }

\begin{document}
\maketitle
\begin{abstract}
A pure pair in a graph $G$ is a pair $(Z_1,Z_2)$ of disjoint sets of vertices such that either every vertex in $Z_1$ is 
adjacent to every vertex in $Z_2$, or there are no edges between $Z_1$ and $Z_2$.  With Maria Chudnovsky, we recently proved 
that, for every forest $F$, 
every graph $G$ with at least two vertices that does not contain $F$ or its complement as an induced subgraph has a pure 
pair $(Z_1,Z_2)$ with $|Z_1|,|Z_2|$ linear in $|G|$.   

Here we investigate what we can say about pure pairs in an {\em ordered} graph $G$, when
we exclude an ordered forest $F$ and its complement as induced subgraphs. Fox showed that there need not be a linear pure pair; but Pach and Tomon showed that if $F$ is a monotone path then there is a pure pair of size $c|G|/\log |G|$.  We generalise this to all ordered forests, at the cost of a slightly worse bound:
we prove that, for every ordered forest $F$,
every ordered graph $G$ with at least two vertices that does not contain $F$ or its complement as an induced 
subgraph has a pure pair of size $|G|^{1-o(1)}$.

%
\end{abstract}

\section{Introduction}

In this paper, all graphs are finite and with no loops or parallel edges, and $|G|$ denotes the number of vertices of $G$.
Two disjoint sets are {\em complete} to each other if every vertex of the first is adjacent to every vertex of the second, and
{\em anticomplete} if there are no edges between them. 
A pair $(Z_1,Z_2)$ of subsets of $V(G)$ is {\em pure} if 
$Z_1$ is either complete or anticomplete to $Z_2$. 
A graph $G$ is {\em $H$-free} if no induced subgraph of $G$ is isomorphic to $H$; and if $\mathcal F$ is a family of graphs then a graph is {\em $\mathcal F$-free} if it is $F$-free for all $F\in\mathcal F$.
We denote the complement graph of $H$ by $\overline{H}$.  A {\em hereditary class} or {\em ideal} of graphs is a class of graphs closed under taking induced subgraphs and under isomorphism.

A class $\mathcal G$ of graphs has the {\em strong Erd\H os-Hajnal property} if
there is some $\epsilon>0$ such that every graph $G\in\mathcal G$ with at least two vertices 
contains a pure pair $(A,B)$ such that $|A|,|B|\ge\epsilon|G|$.
Let us consider the class $\mathcal G$ of graphs defined by excluding a finite set $\mathcal F$ of graphs: by considering sparse random graphs, it is easy to show that if the class of $\mathcal F$-free graphs has the strong Erd\H os-Hajnal property then $\mathcal F$ must contain a forest; and by considering complements, it follows also that $\mathcal F$ must contain the complement of a forest.   In an earlier paper~\cite{pure1}, with Maria Chudnovsky, we proved that this is enough to obtain the strong Erd\H os-Hajnal property:
\begin{thm}\label{oldforestsymm}
For every forest $F$, there exists $\epsilon>0$ such that every graph $G$ with at least two vertices that is both $F$-free and $\overline{F}$-free 
contains a pure pair $(Z_1,Z_2)$ with $|Z_1|,|Z_2|\ge \epsilon |G|$.
\end{thm}

We also proved the stronger result that, for sparse graphs, it is enough to exclude just a forest:
\begin{thm}\label{oldmainthm}
For every forest $F$, there exists $\epsilon>0$ such every $F$-free graph $G$ with $|G|\ge2$ has either
\begin{itemize}
\item a vertex with degree at least $\epsilon |G|$; or
\item an anticomplete pair $(Z_1,Z_2)$ with $|Z_1|,|Z_2|\ge \epsilon |G|$.
\end{itemize}
\end{thm}
Again, considering a sparse random graph shows that this does not hold unless $F$ is a forest.

%

In this paper we will be concerned with ordered graphs.
Let us say an {\em ordered graph} is a graph with a linear order on its vertex set, and if $H$ is an ordered graph, $\overline{H}$
denotes the complement graph with the same vertex order. Every induced subgraph inherits an order on
its vertex set in the natural way: let us say an ordered graph $G$ {\em contains}
an ordered graph $H$ if $H$ is isomorphic to an induced subgraph $H'$ of $G$, where the isomorphism carries the order on
$V(H)$
to the inherited order on $V(H')$, and in this case we call $H'$ a {\em copy} of $H$.
We say an ordered graph is {\em $H$-free} if it does not contain the ordered graph $H$.

One could ask for an analogue of \ref{oldforestsymm} for ordered graphs, but it is false. That is a consequence of the following result of 
Fox~\cite{fox}:
\begin{thm}\label{fox}
Let $H$ be the ordered graph with three vertices $h_1,h_2,h_3$ in this order, and with edges $h_1h_2$ and $h_2h_3$.
For all sufficiently large $n$, there is an $H$-free ordered graph $G$ with $n$ vertices, 
such that there is no pure pair $(Z_1,Z_2)$ in $G$ with $|Z_1|,|Z_2|\ge n/\log(n)$.
\end{thm}
To deduce that \ref{oldforestsymm} does not extend to ordered graphs, let $T$ be an ordered forest such that both $T$ and $\overline{T}$
contain $H$; for instance the ordered forest with four vertices $h_1,h_2, h_3, h_4$ in this order, in which $h_1h_2$ and $h_2h_4$ are edges.
Then the graph $G$ of \ref{fox} contains neither $T$ nor its complement.

On the positive side, 
Pach and Tomon \cite{pt1} proved an analogue of \ref{oldmainthm} for monotone paths.
A {\em monotone path} is a path  $x_1\cdots x_k$ with vertices ordered $x_1\le \cdots\le x_k$ (i.e.~the path order agrees with the ordering of the graph).
Pach and Tomon showed that the bound of \ref{fox} is in fact sharp for ordered paths
(see Fox, Pach and T\'oth \cite{fpt} and Fox \cite{fox} for earlier work):

\begin{thm}\label{pttheorem}
Let $P$ be a monotone path.
There exists $\vare>0$ such that every $P$-free ordered graph $G$ with
at least two vertices has either
\begin{itemize}
\item a vertex with degree at least $\vare|G|$; or 
\item an anticomplete pair $(Z_1,Z_2)$ such that $|Z_1|, |Z_2|\ge \vare|G|/\log(|G|)$.
\end{itemize}
\end{thm}

In this paper we prove an analogue of \ref{oldmainthm} that holds for all ordered forests.  We show that excluding any ordered forest guarantees either a vertex of linear degree or an anticomplete pair of size $|G|^{1-o(1)}$.
\begin{thm}\label{sparsethm}
For every ordered forest $T$, and all $c>0$, there exists $\vare>0$ such that every $T$-free ordered graph $G$ with
at least two vertices has either
\begin{itemize}
\item a vertex with degree at least $\vare|G|$; or 
\item an anticomplete pair $(Z_1,Z_2)$ such that $|Z_1|, |Z_2|\ge \vare|G|^{1-c}$.
\end{itemize}
\end{thm}

We had to reduce the $\vare|G|/\log(|G|)$ bound in \ref{pttheorem} to $\vare|G|^{1-c}$ to find a proof, but in fact 
for most ordered trees, this is essentially best possible. Pach and Tomon~\cite{pt2} showed:
\begin{thm}\label{bestposs}
Let $H_1$ be the ordered tree with ordered vertex set $\{v_1\ll v_4\}$ and edges $v_1v_2,v_1v_3,v_1v_4$; and let $H_2$
be the ordered tree with the same ordered vertex set and edges $v_2v_3,v_1v_3,v_1v_4$.
For every $\vare>0$ there exist $\delta > 0$ and $n_0$ with the following property.
For every positive integer $n \ge n_0$, there is an ordered graph $G$ with $n$ vertices and maximum degree at most $\vare n$,
such that
\begin{itemize}
\item  if $A,B\subseteq V(G)$ are anticomplete then $\min(|A|,|B|)\le n^{1-\delta}$; and
\item  $G$ does not contain either of the
ordered trees $H_1,H_2$ as an induced ordered subgraph.
\end{itemize}
\end{thm}

As we will see in the next section, \ref{sparsethm} implies that 
excluding an ordered forest and its complement gives a pure pair of size $|G|^{1-o(1)}$:
\begin{thm}\label{mainthm}
For every ordered forest $T$, and all $c>0$, there exists $\vare>0$ such that 
if $G$ is an ordered graph with $|G|>1$ that is both $T$-free and $\overline{T}$-free, then
$G$ contains a pure pair $(Z_1,Z_2)$ with $|Z_1|, |Z_2|\ge \vare |G|^{1-c}$.
\end{thm}
This characterizes ordered forests and their complements, in that no other ordered graphs $T$ have the property of \ref{mainthm}, because of the following:
\begin{thm}\label{converse}
For every ordered graph $T$ such that neither of $T,\overline{T}$ is a forest, there exists $c>0$ such that for all $\vare>0$,
there are infinitely many ordered graphs $G$
not containing $T$ or its complement, in which there is no pure pair $(Z_1,Z_2)$ with $|Z_1|,|Z_2|\ge \vare|G|^{1-c}$.
\end{thm}
\Proof
Choose an integer $g$ such that both $T$ and $\overline{T}$ have a cycle of length at most $g$. Let $c<1/g$, and let $\vare>0$. If we take a 
random graph $G$ on $n$ vertices where $n$ is sufficiently large, in which every edge is present independently with probability 
$\frac12 n^{-1+1/g}$, then with high probability, there will be a set $X$ of at least $n/2$ vertices in which $G[X]$ has no cycle
of length at most $g$ (and so contains neither of $T,\overline{T}$) and has no pure pair $Z_1,Z_2$ with $|Z_1|,|Z_2|\ge \vare|X|^{1-c}$.~\bbox

\section{Reduction to the sparse case}

The following very useful result is due to V. R\"odl~\cite{rodl}:

\begin{thm}\label{rodl}
For every graph $H$ and all $\vare>0$ there exists $\delta>0$ such that for every $H$-free graph $G$,
there exists $X\subseteq V(G)$ with $|X|\ge \delta|G|$ such that in one of $G[X]$, $\overline{G}[X]$,
every vertex in $X$ has degree less than $\vare|X|$.
\end{thm}

We will show that the same is true when $G,H$ are ordered graphs, because of the following result of R\"odl and Winkler~\cite{RW}:

\begin{thm}\label{orderedornot} For every ordered graph $H$, there exists a graph $H'$
such that, for
every ordering of $V(H')$, the resulting ordered graph contains $H$.
\end{thm}

A version of \ref{rodl} for ordered graphs follows easily:
\begin{thm}\label{orderedrodl}
For every ordered graph $H$ and all $\vare>0$ there exists $\delta>0$ such that for every $H$-free ordered graph $G$,
there exists $X\subseteq V(G)$ with $|X|\ge \delta|G|$ such that in one of $G[X]$, $\overline{G}[X]$,
every vertex in $X$ has degree less than $\vare|X|$.
\end{thm}
\Proof
Choose $H'$ as in \ref{orderedornot}; and choose $\delta$ as in \ref{rodl} with $H$ replaced by $H'$. If $G$ is an $H$-free 
ordered graph, then the underlying unordered graph is $H'$-free, and so the result holds by the choice of $\delta$.~\bbox

\noindent{\bf Proof of \ref{mainthm}, assuming \ref{sparsethm}.\ \ }
Let $T$ be an ordered forest and $c>0$. Let $\vare'$ satisfy \ref{sparsethm} (replacing $\vare$). Now let $\delta$ 
satisfy \ref{orderedrodl} with $H, \vare$ replaced by $T,\vare'$; and let $\vare=\vare'\delta$. We claim that $\vare$ satisfies \ref{mainthm}.
To see this, let $G$ be an ordered graph with
$|G|\ge 2$ that is $T$-free and $\overline{T}$-free. From the choice of $\delta$, 
there exists $X\subseteq V(G)$ with $|X|\ge \delta|G|$ such that in one of $G[X]$, $\overline{G}[X]$,
every vertex in $X$ has degree less than $\vare'|X|$.  Suppose that $|X|=1$; then $\vare|G|\le \delta|G|\le 1$, and any two vertices
of $G$ make a pure pair of singletons sets that satisfy the theorem. So we may assume that $|X|>1$. By taking complements if necessary,
we may assume that every vertex in $X$ has degree in $G[X]$ less than $\vare'|X|$. By \ref{sparsethm} applied to $G[X]$,
there is an anticomplete pair of subsets of $X$, both of cardinality at least 
$$\vare' |X|^{1-c}\ge \vare'\delta^{1-c} |G|^{1-c}\ge \vare'\delta|G|^{1-c}=\vare|G|^{1-c}.$$
This proves \ref{mainthm}.~\bbox

Actually there is a further small strengthening, the following (eliminating the multiplicative constant $\vare$):
\begin{thm}\label{verysparsethm}
For every ordered forest $T$, and all $c>0$, there exists $\vare>0$ such that, if $G$ is a $T$-free ordered graph with
$|G|> 1/\vare$, then either some vertex has degree at least $\vare|G|$, or
there are disjoint $Z_1,Z_2\subseteq V(G)$ such that $|Z_1|, |Z_2|\ge |G|^{1-c}$ and $Z_1$ is anticomplete to $Z_2$.
\end{thm}
\noindent{\bf Proof of \ref{sparsethm}, assuming \ref{verysparsethm}.\ \ }
Let $T$ be an ordered forest and $c>0$. Let $\vare$ be as in \ref{verysparsethm}; we claim that it also satisfies \ref{sparsethm}. 
Let $G$ be a  $T$-free ordered graph with $|G|\ge 2$. If $|G|>1/\vare$, then the result follows from
\ref{verysparsethm}, so we may assume that $|G|\le 1/\vare$. 
Let $u,v\in V(G)$ be distinct. If they are adjacent then $v$ has degree at least $1\le \vare|G|$; and if they are nonadjacent
then $\{u\},\{v\}$ is 
an anticomplete pair
both of cardinality at least $\vare|G|\ge \vare|G|^{1-c}$. In either case the theorem holds.
This proves \ref{sparsethm}.~\bbox

It remains to prove \ref{verysparsethm}, and that occupies the remainder of the paper.

\section{Blockades, and a sketch of the proof}\label{sec:sketch}

In this section we give an idea of how the proof of \ref{verysparsethm} will go. We work by induction on $|T|$; we choose an
appropriate $\vare>0$; and now we are given
an ordered graph $G$ with at least $1/\vare$ vertices, with maximum degree less than $\vare|G|$, and in which there is no
anticomplete pair of sets both of cardinality at least $|G|^{1-c}$. We need to show that $G$ contains $T$ as an ordered
induced subgraph.

Let $v\in V(T)$ have degree one, and let $T'$ be obtained from $T$ by deleting $v$. From the inductive hypothesis,
there is certainly a copy of $T'$ in $G$,
but we need to produce some vertex of $G$ with the right adjacency to the copy of $T'$, and in the right position in the order
of $G$, to supply the missing leaf of $T$. How can we do this?

It would be helpful if we could arrange that the vertices of the copy of $T'$ are nicely spaced out in the order of $G$; and a
convenient way to do this is to prove a stronger theorem, that if we partition $V(G)$ into some large constant number
of intervals (called ``blocks'') and all the blocks are about the same size, then we can find the ordered copy
of $T'$ with all its vertices in different blocks (we call this being ``rainbow'' with respect to the system of blocks, which we call
a ``blockade''). Of course, now we have to find the missing leaf in some block that has not
been used yet, to carry the induction through.

Let us make some definitions precise.
A {\em blockade} in a graph $G$ is a family $(B_i:i\in I)$ of pairwise disjoint nonempty subsets of $V(G)$,
where $I$ is a set of integers. If the graph is ordered, we also require that the sets $B_i$ are intervals numbered
 in order; that is, if $i,j\in I$ and $i<j$, and $u\in B_i$ and $v\in B_j$, then $u$ is earlier than $v$ in the order of the ordered graph $G$.
We call
the sets $B_i\;(i\in I)$ its {\em blocks}, and $|I|$
its {\em length}.
When $I=\{1\ll k\}$ we sometimes write $(B_1\ll B_k)$ for $(B_i:i\in I)$.
It is convenient not to insist that all blocks have the same cardinality, but what matters is that the smallest
block is not too small. If
the smallest block has cardinality $w$, we call $w$ the {\em width} and $\sigma$ the {\em shrinkage} of the blockade, where
$|G|^{1-\sigma}=w$. (If $I=\emptyset$, the width is $|G|$ and shrinkage is $0$.)

If $\mathcal{B}$ is a blockade in a graph $G$, an induced subgraph $H$ of $G$ is {\em $\mathcal{B}$-rainbow} if
every vertex of $H$ belongs to some block of $\mathcal{B}$, and every block of $\mathcal{B}$ contains at most one vertex of $H$.
If $A,B\subseteq V(G)$ are disjoint, we say $A$ {\em covers} $B$ if every vertex of $B$
has a neighbour in $A$.

Let $X,Y$ be disjoint nonempty subsets of $V(G)$. The {\em max-degree from $X$ to $Y$} is defined to be the maximum
over all $v\in X$ of the number of neighbours of $v$ in $Y$.
Let $(B_i:i\in I)$ be a blockade in a graph $G$, and for all distinct $i,j\in I$ let $d_{i,j}$ be the max-degree from $B_i$ to $B_j$. (Note that $d_{i,j}$ can be quite different from $d_{j,i}$.)
Define $d_{i,i}=0$ for all $i\in I$. We call $d_{i,j}\;(i,j\in I)$ the {\em max-degree function} of the blockade.
Let $\lambda$ be the maximum of $d_{i,j}/|B_j|$, over all distinct $i,j\in I$; we call $\lambda$
the {\em linkage} of $\mathcal{B}$. (If $|I|\le 1$, the linkage is $0$.)

Back to the sketch: if we are aiming at a rainbow copy of $T$, what exactly is the statement we want to prove? We might try:
\\
\\
(1) {\em 
(First attempt.) Let $G$ be an ordered graph such that no two subsets of $V(G)$ of cardinality at least $|G|^{1-c}$ are disjoint and
anticomplete, and every vertex has degree less than $\vare |G|$, where $\vare>0$ is some sufficiently small constant.
If $\mathcal{B}$ is a blockade in $G$
of sufficient length and sufficiently small shrinkage, then there is a $\mathcal{B}$-rainbow copy of $T$.}
\\
\\
Unfortunately, to make the induction on $|V(T)|$ work, the sizes of the blocks of $\mathcal{B}$ need to be sublinear
in $|G|$; and then (1)
is not true, because for instance it might be that every block was complete to every other block. We have to restrict ourselves
to blockades where, although the width might be sublinear in $|G|$, each vertex of each block is only adjacent to a small linear
fraction of each other block, that is, the linkage is at most some fixed constant $<1$. Given this, we can omit the
condition that every vertex has degree less than $\vare |G|$.  So our new goal, and this time we will achieve it, is:
\\
\\
(2) {\em
Let $G$ be an ordered graph such that no two subsets of $V(G)$ of cardinality at least $|G|^{1-c}$ are disjoint and
anticomplete.
If $\mathcal{B}$ is a blockade in $G$,
of sufficient length and and sufficiently small shrinkage and linkage,
then there is a $\mathcal{B}$-rainbow copy of $T$.
}

\bigskip
If all we knew was that there is no anticomplete pair both of cardinality linear in $|G|$, we could not say anything about blockades
in which the block sizes are sublinear; but we have a much stronger hypothesis, that there is no anticomplete pair both of size
at least $|G|^{1-c}$, and that allows us to work with blockades with block sizes down to about $|G|^{1-c}$. In particular, we can
afford to shrink the given blockade by factors of $|G|^{\phi}$ provided that $\phi$ is small compared with $c$.

For (2) to be true, the number of blocks obviously has to be at least $|T|$; but in fact 
the number of blocks needs to be
much bigger, some huge function of $T$ (although independent of the size of $G$). And that brings the problem that, while
we can guarantee that there will be a rainbow copy of $T'$, we do not know ahead of time which particular set
of $|T|-1$ blocks will be used for it. So to prove (2) by induction on $|T|$, we are given a blockade $\mathcal{B}$,
and we would be happy if it had
the property that, wherever the copy of $T'$
appears, we can find the missing leaf. This might not be true; but we will show that we can make it true by 
removing some of the blocks and shrinking the others (not too much).

How can we guarantee that a blockade has this property, that we can always fill in a missing leaf? What we will do is
index the blocks by consecutive integers, say $B_1\ll B_{2k+1}$, and find the copy of $T'$ just using the even blocks,
and find the missing leaf in one of the odd blocks. To be sure we can always do this, we will arrange the property that:
\begin{thm}\label{leafy}
{\bf Desirable property:}
For all odd $i$ and even $j$, there is a subset $X\subseteq B_i$ which covers $B_j$, and there are no edges between $X$ and the other even blocks.
\end{thm}
If we could do that, then whichever set of even blocks is used for $T'$, there is an odd block in the right position in the 
order of $G$, and a vertex in that block with the correct adjacency to the copy of $T'$. So that is our goal;
we start with a blockade, and we want to shrink the even blocks, to make \ref{leafy} true, without shrinking the even 
blocks so much that
there need not be a rainbow copy of $T'$ among these shrunken even blocks.

We will shrink to make \ref{leafy} hold in three steps:

{\bf Step 1.} First, we will arrange (by shrinking the sets $B_i$) that for every pair of blocks $B_i, B_j$, most vertices in $B_i$ have about the same number of 
neighbours in $B_j$, about the max-degree from $B_i$ to $B_j$; and that this remains true even if we shrink $B_i, B_j$ by 
constant factors. This is called being ``shrink-resistant''. This can be accomplished as follows. 
For each pair $i,j$, if we can delete a small 
constant fraction of $B_i$
and $B_j$, and drive down the max-degree from $B_i$ into $B_j$ by a factor of a constant power of $|G|$, say $|G|^\phi$,
we should do so, and repeat. If we adjust the numbers correctly (in particular, choosing $\phi>0$ much smaller than our
target $c$) this process must terminate after a bounded number of steps, or we will find an anticomplete pair of sets
both of size at least $|G|^{1-c}$, which is impossible. When it terminates, the blocks have shrunk from the original,
but only by a constant factor, and we now have a shrink-resistant blockade.

{\bf Step 2.} Second, we will remove some of the blocks and shrink the others (still just by constant factors, so that we remain 
shrink-resistant), to arrange that for some $\tau$, the density of edges between 
$B_i,B_j$ is roughly equal to $\tau$ for every pair of blocks $B_i,B_j$. To do this, we can assume that the density between
each pair is not very small (because otherwise we could find an anticomplete pair disallowed by the hypothesis), and we partition
the possible densities into a bounded number of bands; the values inside each band differ at most by a factor of $|G|^\phi$,
where $\phi$ is some appropriate small constant.
Think of the bands as colours; then we can apply Ramsey's theorem, and find a long subsequence of blocks within which all the 
densities belong to the same band. Choose $\tau$ in the band and remove all the blocks not in the subsequence; 
then for every pair of blocks that remain, the density between
them is  at least $\tau|G|^{-\phi}$ and at most $\tau|G|^\phi$. Consequently, since the density between $B_i,B_j$ is about
equal to the max-degree from $B_i$ to $B_j$ divided by $|B_j|$ (because of the shrink-resistance), we also have good control
of the max-degree between each pair of blocks.

{\bf Step 3.} Now we will shrink the blocks some more to arrange that \ref{leafy} holds. We will sketch this part
later on, in section~\ref{sec:leafcovering}. This step involves more severe shrinking: each even block might shrink to a factor of 
about $|G|^{-\phi}$ of its original size; and we must
still be able to say that in the blockade made by the even blocks, there is a rainbow copy of $T'$. If the new block size
is less than $|G|^{1-c}$, 
the hypothesis about anticomplete pairs tells us nothing; so we must be careful that, while the blocks
must shrink by polynomial factors, they are only small polynomial factors, less than $|G|^c$. But we can arrange 
the numbers correctly to make this work.
This is the content of section \ref{sec:leafcovering}.

Almost the entire proof does not involve the ordering of the
various ordered graphs, so until the final section we have written it in terms of unordered graphs.

\section{Shrink-resistance}

Let $\mathcal{B}=(B_i:i\in I)$ be a blockade in a graph $G$, and let $B_i'\subseteq B_i$ for each $i\in I$, all nonempty;
then $(B_i':i\in I)$ is also a blockade, and we call it a {\em contraction} of $\mathcal{B}$.
If $I'\subseteq I$, then $(B_i:i\in I')$ is also a blockade, called a {\em sub-blockade} of $\mathcal{B}$. 

We will prove in this section that if we are given a blockade in a graph $G$, with sufficiently large length and sufficiently
small shrinkage and linkage, then there is a contraction $(B_1\ll B_k)$ of some sub-blockade,  of any prescribed length $k$, still with small (but slightly larger)
shrinkage and linkage,
where all the numbers $d_{i,j}/|B_j|$ are about the same, and for all distinct $i,j$, many of the vertices in $B_i$
have about $d_{i,j}$ neighbours in $B_j$; and these numbers remain about the same, even if we shrink the blocks further by 
constant factors.

Let 
$\mathcal{B}=(B_i:i\in I)$ be a blockade in a graph $G$, with max-degree function $d_{i,j}\;(i,j\in I)$.
The product of the numbers $d_{i,j}$ for all distinct $i,j\in I$ is called the {\em max-degree product} of $\mathcal{B}$.

Let $0<\phi,\mu\le 1$. We say that $\mathcal{B}$ is
{\em $(\phi,\mu)$-shrink-resistant} if for all distinct $h,j\in I$, and 
for all $X\subseteq B_h$ and $Y\subseteq B_j$ with $|X|\ge \mu|B_h|$ and $|Y|\ge \mu|B_j|$,
the max-degree from $X$ to $Y$ is more than $d_{h,j}|G|^{-\phi}$.
We begin with:

\begin{thm}\label{shrinkresistant}
Let $\mathcal{B}=(B_i:i\in I)$ be a
blockade in a graph $G$, and 
let $0<  \phi,\mu\le 1$. Let $\beta=\mu^{1+\frac{1}{\phi}|I|^2}$.   Then either
\begin{itemize}
\item there exist distinct $h,j\in I$, and $B_h'\subseteq B_h$ and $B_j'\subseteq B_j$ with 
$|B_h'|/|B_h|, |B_j'|/|B_j|\ge \beta$, such that $B_h', B_j'$ are anticomplete; or
\item there is a  $(\phi,\mu)$-shrink-resistant contraction $(B_i':i\in I)$ of $\mathcal{B}$, such that 
$|B_i'|\ge \beta |B_i|$
for each $i\in I$.
\end{itemize}
\end{thm}
\Proof
Let $T=\lfloor \frac{1}{\phi}|I|^2 \rfloor$.
Choose an integer $t$ with $0\le t\le T+1$ maximum such that there is
a contraction $\mathcal{B}'=(B_i':i\in I)$ of $\mathcal{B}$ with
\begin{itemize}
\item $|B_i'|\ge \mu^t|B_i|$ for each $i\in I$; and
\item max-degree product at most $|G|^{|I|^2-\phi t}$.
\end{itemize}
(This is possible since we may take $t=0$ and
$\mathcal{B}'=\mathcal{B}$.) Let $d_{h,j}\;(h,j\in I)$ be the max-degree function of $\mathcal{B}'$.
\\
\\
(1) {\em We may assume that $d_{h,j}\ge 1$ for all distinct $h,j\in I$, and so $t\le T$.}
\\
\\
If $d_{h,j}<1$, then $d_{h,j}=0$, since it is
an integer. Thus $B_h', B_j'$ are anticomplete.
Since $t\le T+1$ and hence $\mu^t\ge \mu^{T+1}\ge \beta$, it follows that $|B_h'|/|B_h|, |B_j'|/|B_j|\ge \beta$, and the
first outcome of the theorem holds. Thus we may assume that $d_{h,j}\ge 1$ and similarly $d_{j,h}\ge 1$. Hence the
max-degree product of $\mathcal{B}'$ is at least one, and since it is at most $|G|^{|I|^2-\phi t}$, it follows
that $|I|^2-\phi t\ge 0$, and so $t\le T$.
This proves (1).
\\
\\
(2) {\em $(B_i':i\in I)$ is $(\phi,\mu)$-shrink-resistant.}
\\
\\
Let $h,j\in I$ be distinct, and let $C_h\subseteq B_h'$ and $C_j\subseteq B_j'$, with
$|C_h|\ge \mu|B_h'|$ and $|C_j|\ge \mu |B_j'|$. Let $d$ be the max-degree from $C_h$ to $C_j$. For all $i\in I$
with $i\ne h,j$ let $C_i=B_i'$. 
From the maximality of $t$, and since $t\le T$, it follows that the
max-degree product of $(C_i:i\in I)$ is more than $|G|^{|I|^2-\phi (t+1)}$. 
Since it is at most $d/d_{h,j}$ times the
max-degree product of $(B_i':i\in I)$, and the latter is at most $|G|^{|I|^2-\phi t}$, it follows that
$d/d_{h,j}> |G|^{-\phi}$. This proves (2).

\bigskip

Since 
$|B_i'|\ge \mu^t|B_i|\ge \beta|B_i|$ for each $i\in I$, the second outcome of the theorem holds. This proves
\ref{shrinkresistant}.~\bbox

Let $(B_i:i\in I)$ be a blockade in a graph $G$, 
and let $0< \tau,\phi,\mu\le 1$. 
We say that $\tau$ is a 
{\em $(\phi,\mu)$-band} for $(B_i:i\in I)$ if 
\begin{itemize}
\item for all distinct $h,j\in I$, the max-degree from $B_h$ to $B_j$ is at most $\tau|B_j|$; and
\item for all distinct $h,j\in I$, and all $X\subseteq B_h$ and $Y\subseteq B_j$ with $|X|\ge \mu|B_h|$ and $|Y|\ge \mu|B_j|$,
the max-degree from $X$ to $Y$ is more than $\tau|G|^{-\phi}|B_j|$.
\end{itemize}
We observe that if $\tau$ is a
{\em $(\phi,\mu)$-band} for $(B_i:i\in I)$, then the linkage of $(B_i:i\in I)$ is at most $\tau$.

\begin{thm}\label{ramsey}
Let $k\ge 0$ be an integer, and let $0< \phi, \mu\le 1$ and $\phi\le 1/5$. Then there is an integer $K\ge k$
with the following property. Let $G$ be a graph, 
and let $(B_i:i\in I)$ be a $(\phi,\mu)$-shrink-resistant blockade in $G$, of length at least $K$. 
Assume that $1-\mu\ge |G|^{-\phi}$.
Then there exists $I'\subseteq I$ with $|I'|=k$ such that $(B_i:i\in I')$ 
has a $(5\phi,\mu)$-band.
\end{thm}
\Proof
From Ramsey's theorem, there is an integer $K\ge 1$ such that for every complete graph with vertex set $I$ where $|I|\ge K$, 
and every colouring of its edges with $ \lfloor 1/(2\phi)+2 \rfloor$ colours, there exists $I'\subseteq I$ 
with $|I'|= k$ such that all edges with both ends in $I'$ have the same colour.

Let $(B_i:i\in I)$ be a $(\phi,\mu)$-shrink-resistant blockade in $G$, where $|I|\ge K$,
with max-degree function $d_{i,j}\;(i,j\in I)$. 
\\
\\
(1) {\em For all distinct $h,j\in I$, there is an integer $t\ge 0$ such that
$$|G|^{-2t\phi}< d_{h,j}/|B_j|, d_{j,h}/|B_h|\le |G|^{-2(t-2)\phi}.$$}
\\
\\
It suffices to show that $(d_{h,j}/|B_j|)|G|^{-2\phi}< d_{j,h}/|B_h|$.
At least $(1-\mu)|B_h|$ vertices in $B_h$ have more than $d_{h,j}|G|^{-\phi}$ neighbours in $B_j$, because otherwise
there would be a set $X\subseteq B_h$ with $|X|\ge \mu|B_h|$ such that the max-degree from $X$ to $B_j$ is at most
$d_{h,j}|G|^{-\phi}$, contradicting that $\mathcal{B}$ is $(\phi,\mu)$-shrink-resistant.
So there are more than $(1-\mu)d_{h,j}|B_h|\cdot|G|^{-\phi}$ edges between $B_h$ and $B_j$.
But there are at most $d_{j,h}|B_j|$ such edges; and so 
$(d_{h,j}/|B_j|)(1-\mu)|G|^{-\phi}< d_{j,h}/|B_h|$. Since $1-\mu\ge |G|^{-\phi}$, it follows that
$(d_{h,j}/|B_j|)|G|^{-2\phi}< d_{j,h}/|B_h|$. Since we also have $(d_{j,h}/|B_h|)|G|^{-2\phi}< d_{h,j}/|B_j|$,
this proves (1). 

\bigskip

For all $h,j\in I$ with $h<j$, let $t\ge 0$ be as in (1);
we call $t$ the {\em type} of the pair $(h,j)$.
We claim that for all such $h,j$, the type $t$ of $(h,j)$ satisfies $0< t\le 1/(2\phi)+2$.
Since $|G|^{-2t\phi}< d_{h,j}/|B_j|\le 1$, it follows that $t>0$. 
Since $d_{h,j}\ge 1$ (from the definition of $(\phi,\mu)$-shrink-resistant), and $|B_j|\le |G|$, it follows that
$1/|G|\le d_{h,j}/|B_j|\le |G|^{-2(t-2)\phi}$, and so $1\le  |G|^{1-2(t-2)\phi}$, that is,
$2(t-2)\phi \le 1$. This proves our claim that $0< t\le 1/(2\phi)+2$.
Thus $t$ is one of the
integers $1\ll \lfloor 1/(2\phi)+2\rfloor$.

From the choice of $K$,
there exists $I'\subseteq I$ 
with $|I'|=k$ such that every pair $(h,j)$ with $h<j$ and $h,j\in I'$ has the same type, $t$ say. 
Let $\tau=|G|^{-2(t-2)\phi}$; then for all
distinct $h,j\in I'$,
$$\tau|G|^{-4\phi}\le d_{h,j}/|B_j|\le  \tau.$$
We claim that $\tau$ is a $(5\phi,\mu)$-band for $(B_i:i\in I')$.
To show this, it remains to show that 
for all distinct $h,j\in I'$, and
for all $X\subseteq B_h$ and $Y\subseteq B_j$ with $|X|\ge \mu|B_h|$ and $|Y|\ge \mu|B_j|$,
the max-degree from $X$ to $Y$ is more than $\tau|G|^{-5\phi}|B_j|$.
But $\mathcal{B}$ is $(\phi,\mu)$-shrink-resistant, and so the max-degree from $X$ to $Y$ is more than $d_{h,j}|G|^{-\phi}$;
and since $d_{h,j}\ge \tau|G|^{-4\phi}|B_j|$, the claim follows.
This proves \ref{ramsey}.~\bbox

By combining \ref{shrinkresistant} and \ref{ramsey}, we deduce:

\begin{thm}\label{homog}
Let $k\ge 0$ be an integer, and let  $0<  \phi,\mu\le 1$. Then there exists an integer $K>0$ with the following property.
Let $\mathcal{B}=(B_i:i\in I)$ be a
blockade of length at least $K$ in a graph $G$, where $|G|^{\phi/5}\ge 1/(1-\mu)$.
Let $\beta=\mu^{1+\frac{5}{\phi}|I|^2}$.   Then either
\begin{itemize}
\item there exist distinct $h,j\in I$, and $B_h'\subseteq B_h$ and $B_j'\subseteq B_j$ with
$|B_h'|/|B_h|, |B_j'|/|B_j|\ge \beta$, such that $B_h', B_j'$ are anticomplete; or
\item there exist $I'\subseteq I$ with $|I'|=k$, and a subset $B_i'\subseteq B_i$ for each $i\in I'$, such that $|B_i'|\ge \beta |B_i|$
for each $i\in I'$, and $(B_i':i\in I')$ has a $(\phi,\mu)$-band.
\end{itemize}
\end{thm}
\Proof
Let $K$ satisfy \ref{ramsey} with $\phi$ replaced by $\phi/5$.
Let $G$ be a graph with $|G|^{\phi/5}\ge 1/(1-\mu)$. By \ref{shrinkresistant}, either
\begin{itemize}
\item there exist distinct $h,j\in I$, and $B_h'\subseteq B_h$ and $B_j'\subseteq B_j$ with
$|B_h'|/|B_h|, |B_j'|/|B_j|\ge \beta$, such that $B_h', B_j'$ are anticomplete; or
\item there is a  $(\phi/5,\mu)$-shrink-resistant contraction $\mathcal{B}'=(B_i':i\in I)$ of $\mathcal{B}$, such that
$|B_i'|\ge \beta |B_i|$
for each $i\in I$.
\end{itemize}
In the first case the first outcome of the theorem holds. In the second case,
by \ref{ramsey} applied to $\mathcal{B}'$ with $\phi$ replaced by $\phi/5$, the second outcome of the theorem holds.  This proves \ref{homog}.~\bbox

Consequently we have:
\begin{thm}\label{homog2}
Let $k\ge 0$ be an integer, and let $0<  c,\phi,\mu,\sigma, \Sigma, \Lambda\le 1$ with $\sigma<\Sigma< c$.
Then there exist $\lambda> 0$ and integers $N$ and $K\ge 2$, with the following property.
Let $G$ be a graph with $|G|\ge N$, such that
there do not exist $Z_i\subseteq V(G)$ with $|Z_i|\ge |G|^{1-c}$ for $i = 1,2$, disjoint and anticomplete.
Let $\mathcal{B}=(B_i:i\in I)$ be a
blockade of length at least $K$ in $G$, with shrinkage at most $\sigma$ and linkage at most $\lambda$. Then there exist
$I'\subseteq I$ with $|I'|=k$, and a subset $B_i'\subseteq B_i$ for each $i\in I'$, such that 
$(B_i':i\in I')$ has  shrinkage at most $\Sigma$, and has a $(\phi,\mu)$-band which is at most $\Lambda$.
\end{thm}
\Proof
Let $K$ satisfy \ref{homog}, and let $\beta=\mu^{1+\frac{5}{\phi}K^2}$.
Let $N\ge 0$ such that $N^{\Sigma-\sigma}\ge 1/\beta$, and $N^{\phi/5}\ge 1/(1-\mu)$.
Let $\lambda=\beta \Lambda $.
Let $G$ be a graph with $|G|\ge N$, such that
there do not exist $Z_i\subseteq V(G)$ with $|Z_i|\ge |G|^{1-c}$ for $i = 1,2$, disjoint and anticomplete.
Let $\mathcal{B}=(B_i:i\in I)$ be a
blockade  of length at least $K$ in $G$, with shrinkage at most $\sigma$ and linkage at most $\lambda$. 
If $h,j\in I$ are distinct,
and $B_h'\subseteq B_h$ and $B_j'\subseteq B_j$ with
$|B_h'|/|B_h|, |B_j'|/|B_j|\ge \beta$, then $B_h', B_j'$ are not anticomplete, since
$|B_h'|\ge \beta|B_h|\ge \beta|G|^{1-\sigma}\ge |G|^{1-c}$
and similarly $|B_j'|\ge |G|^{1-c}$. Thus the first outcome of \ref{homog} does not hold.
Since $|G|^{\phi/5}\ge 1/(1-\mu)$, the second outcome of \ref{homog} holds, that is,
there exist $I'\subseteq I$ with $|I'|=k$, and a subset $B_i'\subseteq B_i$ for each $i\in I'$, such that $|B_i'|\ge \beta |B_i|$
for each $i\in I'$, and $(B_i':i\in I')$ has a $(\phi,\mu)$-band. Since $\mathcal{B}$ has shrinkage at most
$\sigma$, and $|B_i'|\ge \beta |B_i|$
for each $i\in I'$, it follows that $(B_i':i\in I')$ has shrinkage at most $\Sigma$, because
$\beta |G|^{1-\sigma} \ge |G|^{1-\Sigma}$.
Also, since $\mathcal{B}$ has linkage at most $\lambda$, and $|B_i'|\ge \beta |B_i|$
for each $i\in I'$, it follows that $(B_i':i\in I')$ has linkage at most $\lambda/\beta=\Lambda$, and therefore there is a 
$(\phi,\mu)$-band for $(B_i':i\in I')$ that is at most $\Lambda$.
This proves \ref{homog2}.~\bbox

\section{Covering with leaves}\label{sec:leafcovering}

Let us continue the sketch of the proof from section \ref{sec:sketch}. We start with some blockade, and we apply \ref{homog2} to it,
and that gives us a contraction of a sub-blockade, still with large length, with linkage and shrinkage only slightly larger
than the original blockade, and with a $(\phi,\mu)$-band $\tau$.
Renumber the blocks by consecutive integers, say $B_1\ll B_{2k+1}$ (we will not remove any more blocks).
We want to arrange that \ref{leafy} holds, but 
let us see first how to arrange that there is a subset $X$ of $B_1$ that covers $B_2$ and has no edges
to the other even blocks (briefly, $(B_1,B_2)$ is ``fixed up''). Here is a method to construct such a set $X$ (for simplicity let us assume that all the blocks
are the same size; in reality some of the numbers involved have to be normalized by dividing by the size of the appropriate
block). Start with the vertex $v_1$ in $B_1$ that has most neighbours in $B_2$; this is about $\tau|B_2|$ (within a factor
of $|G|^\phi$ say). Since all the max-degrees are about the same
(up to the same factor), the number of neighbours of $v_1$ in each even $B_i$ is at most $|G|^\phi$ times its number
of neighbours in $B_2$. Remove its neighbours from $B_2$; by shrink-resistance, there is another vertex $v_2$ in $B_1$
still with about the same number of neighbours in the remainder of $B_2$. And so on, and let $X=\{v_1,v_2,\ldots\}$: we can continue until $B_2$
has been shrunk so much that shrink-resistance is endangered. But stop before that; stop when we have covered and removed
about $|G|^{-\phi}|B_2|/2$ vertices in $B_2$. From the way we chose the sequence $v_1,v_2,\dots$, we know that the
amount of $B_4$ we have covered is only at most $|G|^\phi$ times the amount of $B_2$ we covered, and so is less than half
of $B_4$, and the same for all even blocks.
Now replace $B_2$ by the part of $B_2$ that $X$ covers, and for every other even block, replace it by the part that $X$ does not cover.

This ``fixes up'' the pair $(B_1,B_2)$. We need to fix up similarly all the odd/even pairs of blocks. In this process so far,
the odd blocks have not been changed, and most even blocks have only shrunk to half their size; but $B_2$
has shrunk tremendously, since it was replaced by the portion covered by $X$, which might have been only a fraction
$|G|^{-\phi}$ of the original block. This is a problem, since shrinking a block that much destroys the shrink-resistance.
The blockade consisting of the blocks different from $B_2$ is still shrink-resistant, but not if we include $B_2$ as well. And shrink-resistance
was used crucially to fix up the pair $(B_1,B_2)$; how can we fix up the pair $(B_3,B_2)$ now that shrink-resistance has gone?

The answer is, to fix up all the pairs involving $B_2$ simultaneously; and then start fixing up the pairs involving $B_4$,
and so on. At a general step of this process, there will be some even values (say a set $H$) such that for all odd $i$,
the pair $(B_i,B_h)$ has been fixed up; the other
even values (say $I$), that have not yet been fixed up with the odd blocks; and the odd blocks themselves, $B_j\:(j\in J)$
say.
For $i\in H$, $B_i$ has been shrunk by a factor of something like $|G|^\phi$; for $i\in I$, $B_i$ has only been shrunk by a
constant factor; and for $j\in J$, $B_j$ has not been shrunk at all. Consequently the blockade formed by
the blocks $B_i\;(i\in I\cup J)$ is still shrink-resistant, and the densities between the pairs of its blocks are all still about $\tau$.
In this general step, \ref{matchlemma} below, we will choose some element
$g$ of $I$, and fix up the pair $(B_j,B_g)$ for all odd blocks $B_j$, and shrink $B_g$, and move $g$ from $I$ into $H$.

There is another thing to keep track of: we want to obtain a rainbow copy of $T'$ among the even blocks after contraction, and 
for this we need to make sure that no vertex in one even block has too many neighbours in another. So, even after $g$ has been 
moved from $I$ to $H$ and $B_g$ has been shrunk by the polynomial factor, we still need to keep track of the max-degree between $B_g$ and the other even blocks.

That was an attempt to explain what is happening in the results of this section. There are unfortunately a lot of parameters involved,
and the formal statement of the results is a little fearsome. Our objective in this section is to prove the following:

\begin{thm}\label{fullleafcover}
Let $k\ge 0$ be an integer, and let $0<c, \sigma, \sigma',\lambda'\le 1$ with $\sigma<\sigma'<c$.
Then there exist $\lambda>0$ and integers $K,N>0$ with the following property.
Let $G$ be a graph with $|G|\ge N$, such that
there do not exist $Z_i\subseteq V(G)$ with $|Z_i|\ge |G|^{1-c}$ for $i = 1,2$, disjoint and anticomplete.
Let $\mathcal{A}=(A_i:i\in I)$ be a
blockade of length at least $K$ in $G$, with shrinkage at most $\sigma$ and linkage at most $\lambda$.
Then there exists $I'\subseteq I$ with $|I'|=k$, such that for every partition $(H,J)$ of $I'$,
there exists $B_h\subseteq A_h$ for each $h\in H$,
where
\begin{itemize}
\item $(B_h:h\in H)$ has shrinkage at most $\sigma'$ and linkage at most $\lambda'$; and
\item for all $h\in H$ and
all $j\in J$ there exists $X\subseteq A_j$ that covers $B_h$ and is anticomplete to $B_{i}$ for all $i\in H\setminus \{h\}$.
\end{itemize}
\end{thm}

To prove \ref{fullleafcover} we proceed in several steps.
We begin with: 

\begin{thm}\label{matchlemma}
Let $k\ge 0$ be an integer, and let $0< \tau, \phi,\mu\le 1$, with $2k\mu\le 1$, and $\phi\le 1/2$, and $4k^2\tau\le 1$.
Let $\{0\},H,I,J$ be pairwise disjoint sets of integers, with union of cardinality $k$.
Let $G$ be a graph, such that $|G|^{\phi}\ge 2(16k^2)^k$.
Let $\mathcal{A}=(A_i:i\in \{0\}\cup H\cup I\cup J)$ be a
blockade in $G$, such that:
\begin{itemize}
\item $\tau$ is a $(\phi,\mu)$-band for $(A_i:i\in \{0\}\cup I\cup J)$; and
\item for each $h\in H$, and each $i\in \{0\}\cup I\cup J$, the max-degree from $A_h$ to $A_i$ is at most $\tau|A_i|$.
\end{itemize}
Then for all $i\in \{0\}\cup H\cup I\cup J$ there exists $B_i\subseteq A_i$, such that:
\begin{itemize}
\item $|B_0|\ge |G|^{-k\phi}|A_0|$, and $|B_i|\ge |A_i|/2$ for all $i\in H\cup I$, and $B_i=A_i$ for all $i\in J$;
\item for all $j\in J$ there exists $C_j\subseteq A_j$ that covers $B_{0}$ and is anticomplete to all the sets
$B_i\; (i\in H\cup I)$;
\item $2\tau$ is a $(2\phi,2\mu)$-band for $(B_i:i\in I\cup J)$;
\item for each $h\in H\cup \{0\}$, and each $i\in I\cup J$, the max-degree from $B_h$ to $B_i$ is at most $2\tau|B_i|$; and
\item for each $h\in H$, the max-degree from $B_0$ to $B_h$ is at most $4k\tau|B_h|$.
\end{itemize}
\end{thm}
\Proof We may assume that $k\ge 2$.
\\
\\
(1) {\em There exists $D_j\subseteq A_j$ with $|D_j|\ge |A_j|/2$ for each $j\in J\cup\{0\}$, such that:
\begin{itemize}
\item for each $j\in J\cup \{0\}$ and each $h\in H$, every vertex in
$D_j$ has fewer than $2k\tau |A_h|$ neighbours in $A_h$; and
\item for each $j\in J$, every vertex in $D_0$ has
more than $\tau|G|^{-\phi}|A_j|$ neighbours in $D_j$.
\end{itemize}}
\noindent
For each $j\in J\cup \{0\}$ and each $h\in H$, let $Z_{j,h}$ be the set of vertices in $A_j$ that have 
at least $2k\tau|A_h|$ neighbours in $A_h$; and for each $j\in J\cup \{0\}$, let 
$$D_j'=A_j\setminus \bigcup_{h\in H}Z_{j,h}.$$
Let $D_j=D_j'$ for each $j\in J$ (we will choose $D_0\subseteq D_0'$ later). 
Since every vertex in $A_h$ has at most $\tau |A_j|$ neighbours in $A_j$ (by assumption), there are at most $\tau|A_h|\cdot |A_j|$
edges between $A_h$ and $A_j$, and so $2k\tau|A_h|\cdot |Z_{j,h}|\le \tau|A_h|\cdot |A_j|$, that is,
$|Z_{j,h}|\le |A_j|/(2k)$. For each $j\in J$, the union of the sets $Z_{j,h}$ (over all $h\in H$) has cardinality
at most $|A_j|/2$, and so $|D_j|\ge  |A_j|/2$ for each $j\in J$, and the first statement of (1) holds.

For each $j\in J$, let $Z_j$ be the set of vertices in $A_0$ that have at most $\tau|G|^{-\phi}|A_j|$ neighbours in $D_j$.
Since $\tau$ is a $(\phi,\mu)$-band for $(A_i:i\in I\cup J\cup \{0\})$, and $|D_j|\ge |A_j|/2\ge \mu|A_j|$,
it follows that $|Z_j|\le \mu|A_0|\le |A_0|/(2k)$.
Thus the union of the sets $Z_{0,h}\;(h\in H)$ and the sets $Z_j\;(j\in J)$ has cardinality at most $|A_0|/2$, since 
$|H\cup J|\le k$. Let 
$$D_0=A_j\setminus \left(\bigcup_{h\in H}Z_{0,h}\cup \bigcup_{j\in Z} Z_j\right).$$
Thus $|D_0|\ge |A_0|/2$; and 
the second statement of (1) holds. This proves (1).
\\
\\
(2) {\em Let $Y\subseteq D_0$ and $j\in J$. Then there exists $Y'\subseteq Y$ with $|Y'|\ge |G|^{-\phi}|Y|/(16k^2)$, and 
a subset $C_j\subseteq D_j$, such that $C_j$ covers $Y'$, and for each $h\in H\cup I$, 
at most $|A_h|/(2k)$ vertices in $A_h$ have a neighbour in $C_j$.}
\\
\\
We may assume that $Y\ne \emptyset$. 
Every vertex in $Y$ belongs to $D_0$, and hence has more than $\tau|G|^{-\phi}|A_j|$ neighbours in $D_j$. Choose 
$X\subseteq D_j$ maximal such that 
\begin{itemize}
\item $|X|\le 1/(4k^2\tau)$; and
\item $|Y'|\ge (\tau/2)|G|^{-\phi}|X|\cdot |Y|$, where $Y'$ is the set of vertices in $Y$ that have a neighbour in $X$.
\end{itemize}
For each $i\in I$, since $|X|\le 1/(4k^2\tau)$ and every vertex in $X$ has at most $\tau|A_i|$ 
neighbours in $A_i$, it follows that at most $|A_i|/(4k^2)\le |A_i|/(2k)$ vertices in $A_i$ have a neighbour in $X$. For each $h\in H$,
since every vertex in $D_j$ has at most $2k\tau |A_h|$ neighbours in $A_h$, it follows that at most 
$$2k\tau |A_h|\cdot |X|\le 2k\tau |A_h|/(4k^2\tau)=|A_h|/(2k)$$
vertices in $A_h$ have a neighbour in $X$. 
Thus if $|Y'|\ge |Y|/2$ then (2) holds with $C_j=X$, since $1/2\ge |G|^{-\phi}/(16k^2)$; so we may assume that $|Y'|<|Y|/2$, and 
hence $|Y\setminus Y'|\ge |Y|/2$.
Every vertex in $Y\setminus Y'$ has at least $\tau|G|^{-\phi}|A_j|$ neighbours in $D_j$,
and none of these neighbours is in $X$ since $Y\setminus Y'$ is anticomplete to $X$. Thus there exists $v\in D_j\setminus X$
with at least 
$$\tau|G|^{-\phi}|A_j|\frac{|Y|/2}{|D_j\setminus X|}\ge (\tau|Y|/2)|G|^{-\phi}$$
neighbours in $Y\setminus Y'$ (since $|A_j|\ge |D_j\setminus X|$). From the maximality of $X$, replacing $X$ by $X\cup \{v\}$ violates one of the
two bullets in the definition of $X$.  The second is satisfied, and so the first is violated, that is, 
$|X|+1> 1/(4k^2\tau)$. Consequently $X\ne \emptyset$, and so $2|X|\ge |X|+1>1/(4k^2\tau)$, and therefore $|X|>1/(8k^2\tau)$. Since
$|Y'|\ge (\tau/2)|G|^{-\phi}|X|\cdot |Y|$, it follows that $|Y'|> \tau|G|^{-\phi}|Y|/(16k^2\tau)=|G|^{-\phi}|Y|/(16k^2)$.
This proves (2).

\bigskip

By $|J|$ applications of (2), one for each $j\in J$, applied initially with $Y=D_0$, we obtain
that there exists $B_0\subseteq D_0$ with $|B_0|\ge |G|^{-|J|\phi}(16k^2)^{-|J|}|D_0|$, and for each $j\in J$ there exists a subset 
$C_j\subseteq D_j$, such that $C_j$ covers $B_0$, and for each $h\in H\cup I$,
at most $|A_h|/(2k)$ vertices in $A_h$ have a neighbour in $C_j$.

Since $|J|<k$, and $(16k^2)^{-|J|}\ge 2|G|^{-\phi}$, we obtain $|B_0|\ge 2|G|^{-k\phi}|D_0|\ge |G|^{-k\phi}|A_0|$.
For each $h\in H\cup I$, let $B_h$ be the set of vertices in $A_h$ with no neighbours in any of the sets $C_j\;(j\in J)$.
Then $|B_h|\ge |A_h|/2$. 

The conclusion of \ref{matchlemma} has five bullets,  and we have shown that the first two hold.
For the third bullet, we need the following:
\\
\\
(3) {\em Let $B_i=A_i$ for $i\in J$; then $2\tau$ is a $(2\phi,2\mu)$-band for $(B_i:i\in I\cup J)$, and $(B_i:i\in I\cup J)$
has linkage at most $2\tau$.}
\\
\\
Since $\tau\le 1/(4k^2)\le 1/2$, certainly $2\tau\le 1$.
We claim that $2\tau$ is a $(2\phi,2\mu)$-band for $(B_i:i\in I\cup J)$. To show this, let $i,j\in I\cup J$ be distinct; 
we must show that:
\begin{itemize}
\item the max-degree from $B_i$ to $B_j$ is at most $2\tau|B_j|$; and
\item for all $X\subseteq B_i$ and $Y\subseteq B_j$ with $|X|\ge 2\mu|B_i|$ and $|Y|\ge 2\mu|B_j|$,
the max-degree from $X$ to $Y$ is more than $2\tau|G|^{-2\phi}|B_j|$.
\end{itemize}
Since the max-degree from $A_i$ to $A_j$ is at most $\tau|A_j|$, and $|B_j|\ge |A_j|/2$, it follows that the max-degree 
from $B_i$ to $B_j$ is at most $\tau|A_j|\le 2\tau|B_j|$. Since this holds for all distinct $i,j\in I\cup J$, it follows
that $(B_i:i\in I\cup J)$
has linkage at most $2\tau$. Now let $i,j\in I\cup J$ be distinct, and let $X\subseteq B_i$ and $Y\subseteq B_j$ 
with $|X|\ge 2\mu|B_i|$ and $|Y|\ge 2\mu|B_j|$. Thus $|X|\ge \mu|A_i|$ and $|Y|\ge \mu|A_j|$; and since $\tau$
is a  $(\phi,\mu)$-band for $(A_i:i\in I\cup J)$, it follows that the max-degree from $X$ to $Y$ is more than 
$\tau|G|^{-\phi}|A_j|\ge 2\tau|G|^{-2\phi}|B_j|$ since $|A_j|\ge |B_j|$ and $|G|^{-\phi}\le 1/2$.
This proves that $2\tau$ is a $(2\phi,2\mu)$-band for $(B_i:i\in I\cup J)$.
This proves (3).

\bigskip

Consequently the third bullet of the conclusion of the theorem is satisfied. The fourth holds, since for 
$h\in H\cup \{0\}$, and $i\in I\cup J$, every vertex in $B_h$ has at most $\tau|A_i|\le 2\tau|B_i|$ neighbours in $B_i$.
And the fifth holds since every vertex in $B_0$ belongs to $D_0$, and so has at most 
$2k\tau |A_h|\le 4k\tau|B_h|$
neighbours in $A_h$. This proves \ref{matchlemma}.~\bbox

In order to use \ref{matchlemma} we need the following definition.  
Let $k\ge 0$ be an integer, and let $H,I,J$ be disjoint sets of integers, with union of cardinality $k$. 
Let $0< \tau, \phi,\mu,\lambda\le 1$.
Let $\mathcal{B}=(B_i:i\in H\cup I\cup J)$ be a blockade in a graph $G$.
Suppose that:
\begin{itemize}
\item $(B_h:h\in H)$ has width at least $w$ and linkage at most $\lambda$;
\item $(B_i:i\in I)$ has width at least $W$;
\item for each $h\in H$ and $j\in J$, there exists $X\subseteq B_j$ such that $X$ covers $B_h$ and is anticomplete to $B_i$ for all
$i\in (H\cup I)\setminus \{h\}$;
\item $\tau$ is a $(\phi,\mu)$-band for $(B_i:i\in I\cup J)$;
\item for each $h\in H$, and each $i\in I\cup J$, the max-degree from $B_h$ to $B_i$ is at most $\tau|B_i|$.
\end{itemize}
In these circumstances we say that $\mathcal{B}$ is {\em leaf-covered with partition $(H,I,J)$ and parameters}
$$w,W,\lambda,\phi,\mu,\tau.$$
From \ref{matchlemma} we deduce:

\begin{thm}\label{moreleaves}
Let $k\ge 0$ be an integer, and let $0< \tau, \phi,\lambda, \mu\le 1$, with $2k\mu\le 1$, and $\phi\le 1/2$, and $4k^2\tau\le 1$, 
and $\lambda\ge 2k\tau$.
Let $G$ be a graph, such that $|G|^{\phi}\ge 2(16k^2)^k$.
Let $\mathcal{A}$ be a blockade of length $k$ in $G$ that is leaf-covered with partition $(H,I,J)$ and parameters 
$$w,W,\lambda,\phi,\mu,\tau,$$ 
where $|G|^{-k\phi}W\ge w/4$. Suppose that $g\in I$. Then there is a contraction $\mathcal{B}=(B_i:i\in H\cup I\cup J)$ of 
$\mathcal{A}$, such that $\mathcal{B}$ is leaf-covered with partition $(H\cup \{g\},I\setminus \{g\}, J)$ and
parameters 
$$w/4,W/2,4\lambda,2\phi,2\mu,2\tau.$$
\end{thm}
\Proof
We may assume that $g=0$. Thus $H,I\setminus \{0\},J,\{0\}$ are pairwise disjoint with union of cardinality $k$.
By \ref{matchlemma} with $I$ replaced by $I\setminus \{0\}$, for all $i\in H\cup I\cup J$ there exists $A_i'\subseteq A_i$, such that:
\begin{itemize}
\item $|A_0'|\ge |G|^{-k\phi}|A_0|$, and $|A_i'|\ge |A_i|/2$ for all $i\in (H\cup I)\setminus \{0\}$, and $A_j'=A_j$ for all $j\in J$;
\item for all $j\in J$ there exists $C_j\subseteq A_j$ that covers $A_{0}'$ and is anticomplete to all the sets
$A_i'\; (i\in (H\cup I)\setminus \{0\})$;
\item $2\tau$ is a $(2\phi,2\mu)$-band for $(A_i':i\in (I\cup J)\setminus \{0\})$;
\item for each $h\in H\cup \{0\}$, and each $i\in (I\cup J)\setminus \{0\}$, the max-degree from $A_h'$ to  $A_i'$ is at most
$2\tau|A_i'|$; and
\item for each $h\in H$, the max-degree from $A_0'$ to $A_h'$ is at most $4k\tau|A_h'|$.
\end{itemize}
Let $B_i=A_i'$ for each $i\in I\cup J$.
For each $h\in H$, let $B_h$ be the set of vertices in $A_h'$ that have at most $8k\tau|B_0|$ neighbours in $B_0$. 
Since there are at most $4k\tau|A_h'|\cdot |B_0|$ edges between $A_h'$ and $B_0$, it follows that $|B_h|\ge |A_h'|/2$.
We claim that $(B_i:i\in H\cup I\cup J)$ is leaf-covered with partition $(H\cup \{0\},I\setminus \{0\}, J)$ and
parameters 
$$w/4,W/2,4\lambda,2\phi,2\mu,2\tau.$$ 
To show this, we must check that:
\begin{itemize}
\item $(B_h:h\in H\cup \{0\})$ has width at least $w/4$ and linkage at most $4\lambda$;
\item $(B_i:i\in I\setminus \{0\})$ has width at least $W/2$;
\item for each $h\in H\cup \{0\}$ and $j\in J$, there exists $X\subseteq B_j$ such that $X$ covers $B_h$ and is anticomplete to $B_i$ for all
$i\in (H\cup I)\setminus \{h\}$;
\item $2\tau$ is a $(2\phi,2\mu)$-band for $(B_i:i\in (I\cup J)\setminus \{0\})$; and
\item for each $h\in H\cup \{0\}$, and each $i\in (I\cup J)\setminus \{0\}$, the max-degree from $B_h$ to $B_i$ is at most
$2\tau|B_i|$.
\end{itemize}
For the first bullet: for $h\in H$, $|B_h|\ge |A_h'|/2\ge |A_h|/4\ge w/4$, and 
$$|B_0|\ge |G|^{-k\phi}|A_0|\ge |G|^{-k\phi}W\ge w/4,$$
so $(B_h:h\in H\cup \{0\})$ has width at least $w/4$. Since $(A_h:h\in H)$ has linkage at most $\lambda$, it follows that
$(B_h:h\in H)$ has linkage at most $4\lambda$ (because $|B_h|\ge |A_h|/4$ for each $h\in H$). For each $h\in H$, 
every vertex in $B_0$ has at most
$4k\tau|A_h'|$ neighbours in $A_h'$, and hence at most $8k\tau|B_h|\le 4\lambda|B_0|$ neighbours in $B_h$; and every vertex in $B_h$
has at most $8k\tau|B_0|\le 4\lambda|B_0|$ neighbours in $B_0$. Thus the linkage of $(B_h:h\in H\cup \{0\})$ is at most $4\lambda$.
This proves that the first bullet holds.

The second bullet holds since $|B_i|\ge |A_i|/2\ge W/2$ for each $i\in I\setminus \{0\}$.
The third bullet holds, since if $h\in H$ the statement is true by hypothesis, and if $h=0$ then the statement is true because
we may take $X=C_j$.
The fourth and fifth bullets holds by the application of \ref{matchlemma}. 
This proves \ref{moreleaves}.~\bbox

By repeatedly moving elements from $I$ to $H$ using \ref{moreleaves}, we deduce:
\begin{thm}\label{leafcover}
Let $k\ge 0$ be an integer, and let $I,J$ be disjoint sets of integers with union of cardinality $k$.
Let $0< \tau, \phi,\mu\le 1$, with $k2^k \mu\le 1$, and $\phi2^k\le 1$, and $k^22^{k+1}\tau\le 1$.
Let $G$ be a graph with $|G|^{\phi}\ge 2(16k^2)^k$.
Let $\mathcal{A}=(A_i:i\in I\cup J)$ be a blockade in $G$ with 
a $(\phi,\mu)$-band $\tau$. Let $W$ be the 
width of $(A_i:i\in I)$.
For all $H\subseteq I$, there
is a contraction $\mathcal{B}=(B_i:i\in I\cup J)$ of
$\mathcal{A}$, such that $\mathcal{B}$ is leaf-covered with partition $(H,I\setminus H, J)$ and
parameters 
$$4^{-|H|}|G|^{-k2^{k-1}\phi}W,\; 2^{-|H|}W,\; 4^{|H|}(2k\tau),\; 2^{|H|}\phi,\; 2^{|H|}\mu,\; 2^{|H|}\tau.$$
\end{thm}
\Proof
We proceed by induction on $|H|$. The result is true when $H=\emptyset$, since $\tau$ is a $(\phi,\mu)$-band for 
$\mathcal{A}=(A_i:i\in I\cup J)$. Thus we assume that $H\ne \emptyset$. Choose $g\in H$. From the inductive hypothesis,
there
is a contraction $\mathcal{B}'=(B_i':i\in I\cup J)$ of
$\mathcal{A}$, such that $\mathcal{B}'$ is leaf-covered with partition $(H\setminus \{g\},I\setminus (H\setminus \{g\}), J)$ 
and
parameters 
$$w'=4^{1-|H|}|G|^{-k2^{k-1}\phi}W,\; W'=2^{1-|H|}W,\; \lambda'=4^{|H|-1}(2k\tau),\; \phi'=2^{|H|-1}\phi,\; \mu'=2^{|H|-1}\mu, \;
 \tau'=2^{|H|-1}\tau.$$
Since 
\begin{eqnarray*}
2k\mu'= 2k2^{|H|-1}\mu&\le& 2k2^{k-1}\mu\le 1\\
\phi'=2^{|H|-1}\phi&\le& 2^{k-1}\phi\le 1/2\\
4k^2\tau'=4k^22^{|H|-1}\tau&\le& 4k^22^{k-1}\tau\le 1\\
\lambda'=4^{|H|-1}(2k\tau)&\ge& 2^{|H|-1}(2k\tau)=2k\tau'\\
|G|^{\phi'}=|G|^{2^{|H|-1}\phi}&\ge& |G|^{\phi}\ge 2(16k^2)^k\\
|G|^{-k\phi'}W'= |G|^{-k2^{|H|-1}\phi} 2^{1-|H|}W&\ge& |G|^{-k2^{k-1}\phi} 4^{-|H|}W =w'/4
\end{eqnarray*}
it follows from \ref{moreleaves} that
there is a contraction $\mathcal{B}=(B_i:i\in I\cup J)$ of
$\mathcal{B}'$, such that $\mathcal{B}$ is leaf-covered with partition $(H,I\setminus H, J)$ and
parameters 
$$w'/4,W'/2,4\lambda',2\phi',2\mu',2\tau'.$$ 
But then $\mathcal{B}$ is the required contraction of $\mathcal{A}$.
This proves \ref{leafcover}.~\bbox

We will only apply \ref{leafcover} when $H=I$, and in that case it becomes much simpler, so much so that it is worth 
stating separately, in the following:
\begin{thm}\label{simpleleafcover}
Let $k\ge 0$ be an integer, and let $H,J$ be disjoint sets of integers with union of cardinality $k$.
Let $0< \tau, \phi,\mu\le 1$, with $k2^k \mu\le 1$, and $\phi2^k\le 1$, and $k^22^{k+1}\tau\le 1$.
Let $G$ be a graph with $|G|^{\phi}\ge 2(16k^2)^k$.
Let $\mathcal{A}=(A_i:i\in H\cup J)$ be a blockade in $G$, with a $(\phi,\mu)$-band $\tau$. Let $W$ be the
width of $(A_i:i\in H)$.
For each $h\in H$ there exists $B_h\subseteq A_h$, such that:
\begin{itemize}
\item $(B_h:h\in H)$ has width at least $4^{-|H|}|G|^{-k2^{k-1}\phi}W$ and linkage at most $4^{|H|}(2k\tau)$; and
\item for all $h\in H$ and
all $j\in J$ there exists $X\subseteq A_j$ that covers $B_h$ and is anticomplete to $B_{i}$ for all $i\in H\setminus \{h\}$.
\end{itemize}
\end{thm}
\Proof
By \ref{leafcover}, taking $I=H$,
there
is a contraction $\mathcal{B}=(B_i:i\in H\cup J)$ of
$\mathcal{A}$, such that $\mathcal{B}$ is leaf-covered with partition $(H,\emptyset, J)$ and
parameters
$$4^{-|H|}|G|^{-k2^{k-1}\phi}W,\; 2^{-|H|}W,\; 4^{|H|}(2k\tau),\; 2^{|H|}\phi,\; 2^{|H|}\mu,\; 2^{|H|}\tau.$$
It follows that $(B_h:h\in H)$ has width at least $4^{-|H|}|G|^{-k2^{k-1}\phi}W$ and linkage at most $4^{|H|}(2k\tau)$ (and the 
other four parameters are
irrelevant). This proves \ref{simpleleafcover}.~\bbox

By combining \ref{homog2} and \ref{simpleleafcover}, we obtain the main result of this section, \ref{fullleafcover}, which we restate:
\begin{thm}\label{fullleafcover2}
Let $k\ge 0$ be an integer, and let $0<c, \sigma, \sigma',\lambda'\le 1$ with $\sigma<\sigma'<c$.
Then there exist $\lambda>0$ and integers $K,N>0$ with the following property.
Let $G$ be a graph with $|G|\ge N$, such that
there do not exist $Z_i\subseteq V(G)$ with $|Z_i|\ge |G|^{1-c}$ for $i = 1,2$, disjoint and anticomplete.
Let $\mathcal{A}=(A_i:i\in I)$ be a
blockade of length at least $K$ in $G$, with shrinkage at most $\sigma$ and linkage at most $\lambda$.
Then there exists $I'\subseteq I$ with $|I'|=k$, such that for every partition $(H,J)$ of $I'$,
there exists $B_h\subseteq A_h$ for each $h\in H$, 
where
\begin{itemize}
\item $(B_h:h\in H)$ has shrinkage at most $\sigma'$ and linkage at most $\lambda'$; and
\item for all $h\in H$ and
all $j\in J$ there exists $X\subseteq A_j$ that covers $B_h$ and is anticomplete to $B_{i}$ for all $i\in H\setminus \{h\}$.
\end{itemize}
\end{thm}
\Proof Choose $\Sigma$ with $\sigma<\Sigma<\sigma'$.
Let $\phi$ satisfy $(k2^{k-1}+1)\phi=\sigma'-\Sigma$. Let $\mu=2^{-k}/k$, and $\Lambda=\lambda' 4^{-k}/(2k)$.
Choose  $\lambda> 0$ and integers $N_1$ and $K\ge 2$ such that \ref{homog2} is satisfied with $N$ replaced by
$N_1$. Choose $N\ge N_1$ such that $N^{\phi}\ge 2(16k^2)^k$. We claim that $\lambda, K, N$ satisfy the theorem.

Let  $G$ be a graph with $|G|\ge N$, such that
there do not exist $Z_i\subseteq V(G)$ with $|Z_i|\ge |G|^{1-c}$ for $i = 1,2$, disjoint and anticomplete.
Let $\mathcal{A}=(A_i:i\in I)$ be a
blockade of length at least $K$ in $G$, with shrinkage at most $\sigma$ and linkage at most $\lambda$.
Since $\phi, \Lambda\le 1$, 
it follows from
\ref{homog2} that 
there exist
$I'\subseteq I$ with $|I'|=k$, and a subset $A_i'\subseteq A_i$ for each $i\in I'$, such that
$\mathcal{A}'=(A_i':i\in I')$ has a $(\phi,\mu)$-band $\tau\le \Lambda$, and has  shrinkage at most $\Sigma$.
Let $W=|G|^{1-\Sigma}$; thus $(A_i':i\in I')$ has width at least $W$. 
Since $\phi2^k\le 1$, and 
$$k^22^{k+1}\tau\le k^22^{k+1}\Lambda\le k^22^{k+1}4^{-k}/(2k)\le  1,$$ 
and $|G|^{\phi}\ge 2(16k^2)^k$, \ref{simpleleafcover} implies that
for every partition $(H,J)$ of $I'$, there
exists $B_h\subseteq A_h'$ for each $h\in H$,
such that:
\begin{itemize}
\item $(B_h:h\in H)$ has width at least $4^{-|H|}|G|^{-k2^{k-1}\phi}W$ and linkage at most $4^{|H|}(2k\tau)\le 4^{k}(2k\Lambda)= \lambda'$; and
\item for all $h\in H$ and
all $j\in J$ there exists $X\subseteq A_j'$ that covers $B_h$ and is anticomplete to $B_{i}$ for all $i\in H\setminus \{h\}$.
\end{itemize}
Since $4^{-|H|}\ge 4^{-k}\ge |G|^{-\phi}$, it follows that
$$4^{-|H|}|G|^{-k2^{k-1}\phi}W\ge |G|^{-(1+k2^{k-1})\phi}|G|^{1-\Sigma}= |G|^{1-\sigma'},$$
and so $(B_h:h\in H)$ has shrinkage at most $\sigma'$.
This proves \ref{fullleafcover}.~\bbox

\section{The proof of the main theorem}

In this section we use \ref{fullleafcover} to prove \ref{verysparsethm}. Let $G$ be an ordered graph, and let $H$ be the unordered graph
obtained from $G$ by omitting the ordering. We recall that a {\em blockade} in $G$ is a blockade $(B_i:i\in I)$ in $H$, such that
for all $i,j\in I$ with $i<j$, every vertex of $B_i$ is earlier than each vertex of $B_j$ in the ordering of $G$. Width, shrinkage
and so on are defined as for blockades in unordered graphs.
We will prove:

\begin{thm}\label{rainbow}
Let $0<c\le 1$. For every ordered tree $T$ and all $\sigma$ with $0< \sigma<c$, there exist $\lambda$ with $0< \lambda\le 1$, 
and integers $K,N\ge 0$, with the following property. Let $G$ be an ordered graph with $|G|\ge N$ such that 
there do not exist disjoint $Z_1,Z_2\subseteq V(G)$, where $|Z_1|, |Z_2|\ge |G|^{1-c}$ and $Z_1$ is anticomplete to $Z_2$.
Let $\mathcal{A}$ be a blockade in $G$ of length $K$, with shrinkage at most $\sigma$ and linkage at most $\lambda$.
Then there is an $\mathcal{A}$-rainbow copy of $T$.
\end{thm}
\Proof
We proceed by induction on $|V(T)|$, and we may assume that $|V(T)|\ge 2$. Let $v$ be a leaf of $T$, and let
$T'=T\setminus \{v\}$. Choose $\sigma'$ with $\sigma<\sigma'<c$.
From the inductive hypothesis, there exist $\lambda', K', N'$ so that \ref{rainbow} holds
with $T, \sigma, \lambda, K, N$ replaced by  $T', \sigma', \lambda', K', N'$ respectively. Define $k=2K'+1$.
Choose $\lambda, K,N$ such that \ref{fullleafcover} is satisfied. We claim that $\lambda, K,N$ satisfy \ref{rainbow}.

Let $G$ be an ordered graph with $|G|\ge N$ such that
there do not exist disjoint $Z_1,Z_2\subseteq V(G)$, where $|Z_1|, |Z_2|\ge |G|^{1-c}$ and $Z_1$ is anticomplete to $Z_2$.
Let $\mathcal{A}=(A_i:i\in I)$ be a blockade in $G$ of length $K$, with shrinkage at most $\sigma$ and linkage at most $\lambda$.
From the choice of $\lambda, K,N$, there exists $I'\subseteq I$ with $|I'|=k$, such that for every partition $(H,J)$ of $I'$,
there exists $B_h\subseteq A_h$ for each $h\in H$,
where
\begin{itemize}
\item $(B_h:h\in H)$ has shrinkage at most $\sigma'$ and linkage at most $\lambda'$; and
\item for all $h\in H$ and
all $j\in J$ there exists $X\subseteq A_j$ that covers $B_h$ and is anticomplete to $B_{i}$ for all $i\in H\setminus \{h\}$.
\end{itemize}
Let $I'=\{i_1\ll i_k\}$ where $i_1<i_2<\cdots<i_k$; let $H=\{i_2,i_4,i_6\ll i_{2K'}\}$
and $J=\{i_1,i_3,i_5\ll i_{2K'+1}\}$ (this is well-defined since $k=2K'+1$), and choose $B_h\subseteq A_h$ for each $h\in H$, 
satisfying the two bullets above.
Let $\mathcal{B}=(B_h:h\in H)$.
It follows from the choice of $N', K', \sigma', \lambda'$ that there is a $\mathcal{B}$-rainbow copy of $T'$, and
to simplify notation, we assume that this $\mathcal{B}$-rainbow copy of $T'$ is $T'$ itself.
We recall that $v$ is a leaf of $T$, and $T'=T\setminus \{v\}$. 
Let the linear order of the vertices of $T$ be $(v_1\ll v_n)$, where $v=v_t$. 
Since the vertices of $T'$ appear in the blocks of $\mathcal{B}$ in the correct order, and $J$ interleaves $H$,
there exists
$j\in J$ such that for all $i\in \{1\ll n\}\setminus \{t\}$:
\begin{itemize}
\item if $i<t$ then $v_i\in B_h$ for some $h\in H$ with $h<j$;
\item if $i>t$ then $v_i\in B_h$ for some $h\in H$ with $h>j$.
\end{itemize}
Let $u$ be the neighbour of $v$ in $T$, and let $u\in B_h$.
From the first set of bullets in this proof,
there exists $X\subseteq A_j$ that covers $B_h$ and is anticomplete to $B_{i}$ for all $i\in H\setminus \{h\}$. Choose $v'\in X$
adjacent to $u$. Then adding $v'$ to $T'$ gives a $\mathcal{B}$-rainbow, and hence $\mathcal{A}$-rainbow, copy of $T$.
This proves \ref{rainbow}.~\bbox

Finally we deduce \ref{verysparsethm}, which we restate (and which we have already shown to imply \ref{mainthm}):
\begin{thm}\label{verysparsethm2}
For every ordered forest $T$, and all $c>0$, there exists $\vare>0$ such that, if $G$ is an ordered graph with
$|G|> 1/\vare$, and every vertex has degree less than $\vare|G|$, and there do not exist disjoint
anticomplete sets $Z_1,Z_2$ with $|Z_1|,|Z_2|\ge |G|^{1-c}$, then $G$ contains $T$.
\end{thm}
\Proof
By adding vertices and edges to $T$ if necessary, we may assume that $T$ is an ordered tree.
Let $\sigma=c/2$, and let $\lambda, K, N$ satisfy \ref{rainbow}. 
Choose $M\ge \max(N,K)$ such that $M^{\sigma}\ge 2K$, and let $\vare=\min\left(1/M, \lambda/(2K)\right)$.
Let $G$ be an ordered graph with $|G|> 1/\vare$, such that every vertex has degree less than $\vare|G|$, and there do not exist disjoint
anticomplete sets $Z_1,Z_2$ with $|Z_1|,|Z_2|\ge |G|^{1-c}$.

Since $|G|\ge K$
there is a blockade $\mathcal{B}$ in $G$ of length $K$ and width $W\ge \lfloor |G|/K\rfloor\ge |G|/(2K)$.
Hence $W\ge |G|^{1-\sigma}$, because $|G|\ge M$ and so $ |G|/(2K)\ge |G|^{1-\sigma}$; and therefore $\mathcal{B}$ has shrinkage at most $\sigma$.
Since every vertex has degree less than $\vare|G|$, it follows that $\mathcal{B}$ has linkage at most $\vare|G|/W\le 2K\vare\le \lambda$.
But then from \ref{rainbow} there is a copy of $T$ in $G$. This proves \ref{verysparsethm2} and hence proves \ref{verysparsethm}.~\bbox

\section{Conclusion}

Let us say that a class $\mathcal G$ of graphs or ordered graphs has the {\em Erd\H os-Hajnal property} if
there exists $c>0$ such that every $G\in\mathcal G$ satisfies
$\alpha(G)\omega(G)\ge |G|^c$.   
The Erd\H{o}s-Hajnal conjecture~\cite{EH0, EH} asserts that, for every $H$, the class of $H$-free graphs has the Erd\H os-Hajnal property:
\begin{thm}\label{EHconj}
{\bf Conjecture:} For every graph $H$, there exists $c>0$ such that every $H$-free graph $G$ satisfies
$$\alpha(G)\omega(G)\ge |G|^c.$$
\end{thm}

Alon, Pach and Solymosi \cite{aps} showed that the Erd\H os-Hajnal conjecture is equivalent to the following statement for ordered graphs.

\begin{thm}\label{EHconjordered}
{\bf Conjecture:} For every ordered graph $H$, there exists $c>0$ such that every $H$-free ordered graph $G$ satisfies
$$\alpha(G)\omega(G)\ge |G|^c.$$
\end{thm}

The Erd\H{o}s-Hajnal conjecture (\ref{EHconj}; or equivalently \ref{EHconjordered}) has only been proved for a very small family of graphs.   For example, it remains open for most forests; indeed, it is open even for the five-vertex path.  
However, \ref{oldforestsymm} allows us to say something if we exclude both a forest and its complement.
For an ideal $\mathcal G$ of graphs, the strong Erd\H{o}s-Hajnal property implies the Erd\H{o}s-Hajnal property 
(see~\cite{APPRS,fp}); and the same follows straightforwardly for ideals of ordered graphs.
Thus \ref{oldforestsymm} implies the following:

\begin{thm}\label{oldforestsymmconsequence}
For every forest $F$, the class of graphs that are both $F$-free and $\overline{F}$-free has the Erd\H os-Hajnal property. 
\end{thm}

For {\em ordered} forests, however, the situatation is different: \ref{pttheorem} and \ref{sparsethm} are not strong enough to deduce the Erd\H os-Hajnal property; and we know from \ref{fox} that excluding an ordered forest and its complement is not in general sufficient to obtain the strong Erd\H os-Hajnal property.  Nevertheless, Pach and Tomon \cite{pt2} recently showed the following:

\begin{thm}\label{pathpachtomon}
Let $P$ be a monotone path.  The class of ordered graphs that are both $P$-free and $\overline{P}$-free has the Erd\H os-Hajnal property. 
\end{thm}

It would be interesting to extend this to other ordered trees.  For example, what about extending \ref{pathpachtomon} to all ordered paths?

\end{document}